\documentclass[11pt,leqno]{amsart}
\usepackage{amsfonts,amssymb}
\usepackage{graphicx,tikz}
\usepackage{tikz-cd} 
\usepackage{dsfont}
\usepackage{float}
\usepackage{amsmath, amsthm, amsfonts}
\usepackage{hyperref}
\usepackage{enumitem}  
\usepackage[all,cmtip]{xy}
\usepackage[capitalise]{cleveref}
\usepackage{comment}
\usepackage{enumitem}
\usepackage{amsrefs}
\usepackage{nicefrac}
\usepackage[official]{eurosym}
\usepackage{graphicx}
\usepackage{nomencl}
\usepackage{amsfonts}
\usepackage{hyperref}
\usepackage{tabto}
\usepackage{amssymb}
\usepackage{fancyhdr}
\usepackage{amscd}
\usepackage[english]{babel}

\def\GG{G_{\omega}}

\usepackage[margin=1.15in]{geometry}
\theoremstyle{plain}

\newtheorem*{theorem*}{Theorem}
\newtheorem{thm}{Theorem}[section]

\theoremstyle{definition}

\newtheorem{problem}[thm]{\scshape{Problem}}

\newcommand\T{{\mathcal{T}}}
\makeatletter
\def\cleardoublepage{\clearpage\if@twoside \ifodd\c@page\else
	\hbox{}
	\thispagestyle{empty}
	\newpage
	\if@twocolumn\hbox{}\newpage\fi\fi\fi}
\makeatother

\def\vlongrightarrow{\relbar\joinrel\longrightarrow}

\newcommand{\mapright}[1]{\smash{\stackrel{\text{\tiny{$#1$}}}{\vlongrightarrow}}}

\numberwithin{equation}{section}


\title{Applications of Automaton Groups in Cryptography}
\author[D. Kahrobaei]{Delaram Kahrobaei}
\address{Delaram Kahrobaei: The City University of New York, Departments of Mathematics and Computer Science, Queens College and Initiative for the Theoretical Sciences, CUNY Graduate Center, U.S.A., University of York, U.K. Department of Computer Science, New York University, Tandon School of Engineering, Department of Computer Science and Engineering}
\email{dkahrobaei@qc.cuny.edu}
\author[M. Noce]{Marialaura Noce}
\address{Marialaura Noce: Dipartimento di Matematica, Università degli Studi di Salerno, Italy}
\email{mnoce@unisa.it}

\author[E. Rodaro]{Emanuele Rodaro}
\address{Emanuele Rodaro: Dipartimento di Matematica, Università Politecnico Di Milano, Italy}
\email{emanuele.rodaro@polimi.it}

\begin{document}

\maketitle

\begin{abstract}
In 1991 the first public key protocol involving automaton groups has been proposed. In this paper we give a survey about algorithmic problems around automaton groups which may have potential applications in cryptography. We then present a new public key protocol based on the conjugacy search problem in some families of automaton groups. At the end we offer open problems that could be of interest of group theorists and computer scientists in this direction.
\end{abstract}

\section{Introduction}

 Nowadays the launch of quantum computing has created new challenges in the field of cryptography. The realm of group theory, and more precisely of non-abelian groups offers a collection of examples of new protocols which some of them are claimed to be quantum secure. The field of Group-based Cryptography, is relatively new area in post-quantum cryptography, see the recent survey and book by Kahrobaei et al \cite{kahrobaeinotices, Kahrobaei-BattarbeeBook} for more information.
 
 In cryptography, protocols such as RSA, Diffie-Hellman, and elliptic curve methods depend on the structure of commutative groups and they are related to the difficulty to solve integers factorization and discrete logarithm problem. In 1994, Peter Shor provided an efficient quantum algorithm that solves these problems in polynomial time \cite{Shor}. Since then, researchers were motivated to find alternative methods  to construct new protocols and cryptosystems. 


In this paper, we review what has been done in the direction of using automaton groups as platform for cryptography, and discuss whether they could serve as a post-quantum primitives based on the difficulty of their algorithmic problems.

The class of automaton groups plays an important role in group theory and contains  remarkable examples of infinite groups. These groups provide several solutions to important problems in group theory. As an example, the Grigorchuk group, constructed by Grigorchuk himself in the 1980s, is a counterexample of the General Burnside Problem and is the first example of a group of intermediate word growth, giving a solution to the long-standing problem posted by Milnor in 1960, which asks whether there exists a group that has word growth greater than polynomial but smaller than exponential. 

Automaton groups have already been considered as platforms for some cryptographic protocols. 
The first protocol is based on Grigorchuk groups and it has been introduced in \cite{grigorchukcrypto}, although it was shown to be insecure in \cite{grigorchukcripto2}. In  \cite{grigorchuk2019keyagreement}, Grigorchuk and Grigoriev suggest some specific families of automaton groups, known as branch groups, as platforms for Anshel-Anshel-Goldfeld key-agreement. Among others, they suggest the aforementioned Grigorchuk group and the Basilica group.



This paper is organized as follows. In Section~2 we recall some definitions of automaton groups and basic properties that make these groups suitable for cryptographic applications. In Section~3 we briefly describe the cryptosystems built around automaton groups, namely a public key cryptosystem based on the word decision problem and a key agreement based on the conjugacy search problem. Section~4 is devoted to present a new public key metascheme based on generic automaton groups. In section~5, we discuss the status of the Hidden Subgroup Problem (HSP) for automaton groups, and give some ideas of post-quantum analysis. Finally in Section~6, we conclude the paper by giving open problems for researchers interested to pursue this path of research.

\section{Preliminaries}\label{sec: preliminaries}

\subsection{Automaton groups}
In this section we present the class of automaton groups. For a complete overview of the topic, see \cite{Bartholdi_2010,Grigorchuk_2000,azuk}.

An \emph{alphabet} is a finite set $X$. For each $n\geq 1$, let $X^n$ (resp. $X^{\geq n}$) denote the set of words of length $n$ (resp. of length greater or equal to $n$) over the alphabet $X$ and set $X^0  = \{\emptyset\}$, where $\emptyset$ is the empty word. Moreover we denote by $X^\ast= \bigcup_{n=0}^\infty X^n$, that is the set of all finite words over the alphabet $X$. A Mealy machine or \textit{automaton} is a quadruple $\mathcal{A} =(Q,X,\cdot,\circ)$, where:
\begin{enumerate}
\item $Q$ is a finite set, called the set of states;
\item $X$ is a finite set, called the alphabet;
\item $\lambda: Q\times X \rightarrow Q$ is the restriction map;
\item $\mu: A\times X \rightarrow X$ is the action map.
\end{enumerate}

The automaton $\mathcal{A}$ is \textit{invertible} if, for all $s\in Q$, the transformation $s\circ:=\mu(s,\,):X\rightarrow X$ is a permutation of $X$. An automaton $\mathcal{A}$ can be visually
represented by its \textit{Moore diagram}, that is a directed labeled graph whose vertices are identified with the states of
$\mathcal{A}$. For every state $s\in Q$ and every letter $x\in X$, the diagram has an arrow from $s$ to $s\cdot x:=\lambda(s,x)$
labeled by $x\mid s\circ x$. One can visualize a complete invertible automaton by using such a directed graph: for any $s\in Q$ and $x\in X$ there is exactly one transition of the form 
$$
s\mapright{x|y}t,
$$
for some $t:=s\cdot x\in A$ and $y:=s\circ x \in X$.  
Moreover one can compose the action of the states in $Q$ extending the maps $\circ$ and $\cdot$ to the set $Q^{\ast}$. More precisely, given $w=s_1\cdots s_n\in Q^*$ and $u\in X^{\ast}$ we have
\[
w\circ u=s_2\ldots s_n\circ (s_1\circ u),\quad w\cdot u=\left(s_2\ldots s_n\cdot (s_1\circ u) \right)(s_1\cdot u).
\]
The semigroup generated by this action is called the automaton semigroup defined by $\mathcal A$ and it is denoted by $S(\mathcal{A})$. Clearly, $S(\mathcal{A})$ acts faithfully on $X^{\star}$ via the action $\circ$. Since $\mathcal A$ is invertible we have that for each $s\in Q^{\star}$, $s\circ $ is a bijection, so we may extend $\circ$ to the set of words $(Q\cup Q^{-1})^{\star}$ on the symmetric alphabet $Q\cup Q^{-1}$. In this case we may consider the group generated by $\circ$ that is called automaton group $G(\mathcal A)$ generated by $\mathcal A$. Operationally, the action of each inverse $q^{-1}$, $q\in Q$ is obtained by considering the automaton $\mathcal A^{-1}$ obtained by $\mathcal A$ by swapping input with output: we have a transition $q\mapright{x\mid y} t$ in $\mathcal{A}$ if and only if we have a transition $q^{-1}\mapright{y\mid x} t^{-1}$ in $\mathcal A^{-1}$.
\\
We now give a crucial definition that we will use in our protocol. Suppose that $G$ is a group presented by $\langle Q\mid \mathcal{R}\rangle$. For a word $w\in (Q\cup Q^{-1})^*$, we let $Fact(w)$ denote the set of all factors of $w$, including the empty word $1$. For any $u\in Fact(w)$ there are words $v_1, v_2 \in (Q\cup Q^{-1})^*$ with $w=v_1uv_2$. In this case, we write  $w[u]=v_1^{-1}v_2^{-1}$. Clearly, if $w$ represents the identity in $G$, i.e.,  $w= \mathds{1}$ in $G$, then for any $u\in Fact(w)$ we have $u=w[u]$ in $G$. Given a set $R=\{r_1,\ldots, r_k\}$ of elements representing the identity, we may define a rewriting systems (or semi-Thue system), $\rightarrow_R$ on $(Q\cup Q^{-1})^*$ defined by the binary relation:
\begin{eqnarray*}
& u\rightarrow_R r_i[u], \ \mbox{ for all } u\in Fact(r_i), i=1,\ldots, k;\\
& vv^{-1}\rightarrow_R 1,\ \mbox{ for all } v\in (Q\cup Q^{-1})^*.
\end{eqnarray*}
Note that for any production $x\rightarrow^*_R y$ we have $x=y$ in $G$. In particular, it is easy to see that the equivalence relation $\leftrightarrow^*_R$ generated by $R=\mathcal{R}$ is the word problem (see next subsection for the formal definition) of the group with presentation $\langle Q\mid \mathcal{R}\rangle$.

\subsection{Some algorithmic problems in group theory}

Let $G$ be a group given by a presentation $\langle Q|R \rangle$ where we understand that when we speak of elements of $G$ these are given as a product of generators in $Q^{\pm 1}$. The following three decision problems were introduced by Dehn in 1911. They are defined as follows.

\begin{itemize}
    \item \emph{Word Decision Problem}: For any $g \in G$ written as a product of generators and inverses, determine if there exists an algorithm that verifies whether $g$ is the identity element of $G$.
    \item \emph{Conjugacy Decision Problem}: For any $x, y \in G$, written, as before as a product of generators and inverses, determine if there exists an alogirthm that certifies that $x$ and $y$ are conjugate, which means that there exists an element $z \in G$ such that $x^z=z^{-1}xz=y$.
    \item \emph{Isomorphism Decision Problem}: Let $H$ and $K$ be groups given by finite presentations, determine if there exists an isomorphism from $H$ to $K$. 
\end{itemize}

In general, decision problems involve situations where we need to determine whether a given object, denoted as $\mathcal{O}$, possesses a particular property, denoted as $\mathcal{P}$. On the other hand, for search problems, given
a property $\mathcal{P}$ and an object $\mathcal{O}$, one has to find a witness that shows that $\mathcal{O}$ has the property $\mathcal{P}$. 
Several protocols of non-commutative cryptography, for example Ko-Lee non-commutative Diffie-Hellman protocol \cite{braid}, Anshel-Anshel-Goldfeld Commutator key exchange protocol \cite{AAG}, Kahrobaei-Koupparis digital signature scheme \cite{KK12} and Kahrobaei-Khan non-commutative El Gamal public key encryption scheme \cite{KK06}
are based on difficulty of the conjugacy search problem in the certain proposed groups. For this reason below we present, the ``search'' version of the problems mentioned above. 

As before, we let $G$ be a group given by the presentation $\langle Q|R \rangle$ and we agree that when we refer to an element of $G$ this is given as a product of generators in $Q^{\pm 1}$.

\begin{itemize}
    \item The \emph{Word Search Problem} is: given a finitely presented group $G$ and a word $w =_G 1$  in $G$ find a representation of $w$ as a product of conjugates of defining relators and their inverses.

\item The \emph{Conjugacy Search Problem} is: given a finitely presented group $G$ and $a, b \in G$ such that $a$ and $b$ conjugate, find an element $c \in G$ such that $c^{-1}ac = b$.
\end{itemize}

\subsection{Some examples of automaton groups}

\subsubsection{Grigorchuk automata}

The family of groups $G_{\omega}$ acts on the binary tree $\T_2$ and is defined via the parameter $\omega=(\omega_1, \omega_2, \dots)$, such that  $\omega_i \in \{0,1,2\}$ for any $i\in\mathbb{N}$. Many properties of~$\GG$ depend on the choice of the sequence~$\omega$. For example, the group~$\GG$ is periodic if and only if each $\omega_i \in \{0,1,2\}$ occurs infinitely many times. 
For the sake of simplicity, we only define the most studied one, known as well as the ``first'' Grigorchuk group, which is the group corresponding to the infinite sequence  $\omega=(0,1,2, 0, 1, 2, \dots)$. We write $G$ for $G_{(0,1,2,\dots)}$. Following the notation of automata above, we let the alphabet be $X=\{0,1\}$, the set of states $Q=\{a, b, c, d, e\}$, and the transitions written in graph form (with the state written inside nodes) as follows:
\begin{figure}[h!]
    \centering
    \includegraphics[scale=1]{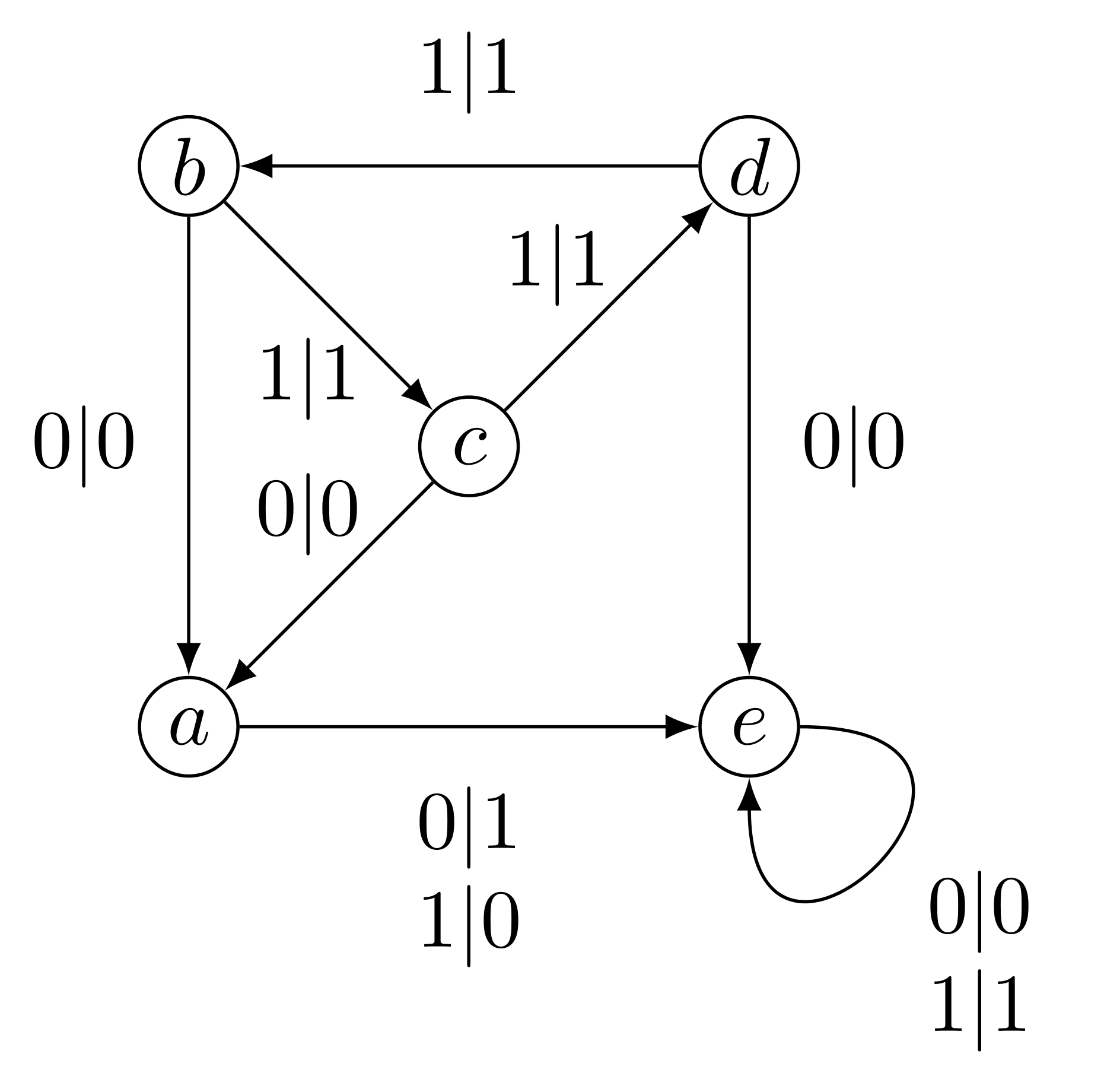}
    \label{fig:grigoautomata}
\end{figure}

\noindent The group $G=\langle a,b,c,d\rangle$ is the (first) Grigorchuk group and possesses a wealth of interesting properties. Among others, it is an infinite finitely generated torsion group providing a counterexample to the General Burnside Problem, and it is a group of intermediate word growth answering a long standing question posed by Milnor.

\subsubsection{Basilica automata}
The Basilica automata $\mathcal{B}$ was introduced by Grigorchuk and Zuk in 2002 \cite{ZUK}, it acts on the binary alphabet and its set of states is $Q=\{e, a, b\}$. Recently, this group has been generalized to the wider class of $p$-Basilica groups \cite{pbasilica}. Below, you can see the $3$-Basilica automaton and what is the action of $a$ and $b$ and the identity $e$ over the alphabet $X=\{0,1,2\}$. Among other properties, we would like to underline that $p$-Basilica groups are torsion-free for any $p$.

\begin{figure}[h!]
    \centering
    \includegraphics[scale=0.8]{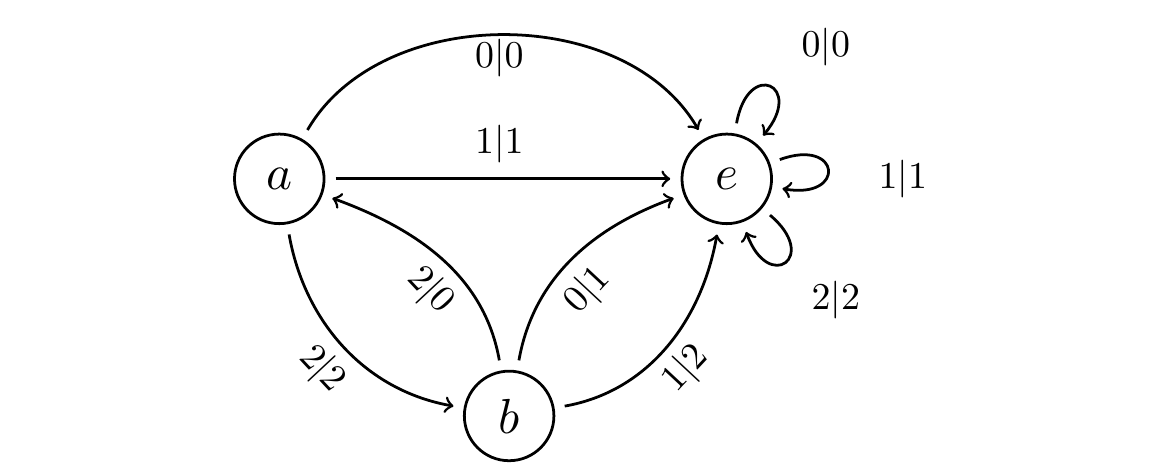}
    \label{fig:basilicaautomata}
\end{figure}

Throughout the paper, for the sake of simplicity, we will refer to Grigorchuk automaton group or Basilica automaton group just as Grigorchuk group and Basilica group.

\subsection{Decidability and complexity}

In 1984, Grigorchuk proved that the Word Decision Problem is solvable in some classes of automaton groups.  Even though it has been shown that many automaton groups have decidable Word Problem, there exist some spinal group that have undecidable Word Problem \cite{Grig1}. 
In 1997 Wilson and Zaleskii solved the Conjugacy Decision Problem for another class of automaton groups, known as GGS $p$-groups for $p$ odd \cite{conjseparability}, while in 2000 Grigorchuk and Wilson proved that the Conjugacy Decition Problem is solvable for a wider class of the so-called branch groups \cite{grigwils}.  Leonov \cite{leonovconj} and Rozhkov \cite{rozhkov}, independently, solved the case $p = 2$. 

Regarding complexity, it has been determined  that the Conjugacy Decision Problem in the Grigorchuk group has polynomial time complexity \cite{miasnykovlys}. In 2017 it has been proved that the Conjugacy Decision Problem in the Grigorchuk group has log-space complexity \cite{log-space}.


We finish this section, by mentioning that the groups presented above, and many examples of known automaton groups are not finitely presented. Neverthless, they admit a recursive presentation, called \emph{endomorphic presentation}. 

For an alphabet~$S$, we denote by $F_S$ the free group on~$S$. A group $G$ has an \emph{$L$-presentation}, or an \emph{endomorphic presentation}, if there exists an alphabet $S$, sets $Q$ and $R$ of reduced words in $F_S$, and a set~$\Phi$ of group homomorphisms $\phi:F_S \rightarrow F_S$ such that $\mathcal{G}$ is isomorphic to a group with the following presentation
\[
\big\langle S \mid  Q\cup \bigcup_{\phi\in \Phi^*} \phi(R) \big\rangle,
\]
where $\Phi^*$ is the monoid generated by~$\Phi$; that is, the closure of $\{1\}\cup \Phi$ under composition.

An $L$-presentation is \emph{finite} if $S$, $Q$, $R$ are finite and $\Phi=\{\phi\}$ consists of just one element.


\section{Cryptographic protocls based on automaton groups}

In this section, we briefly present some cryptosystems known mainly based on the word search and conjugacy search problem in the Grigorchuk groups and in the Basilica groups.

\subsection{Public key cryptosystem  based on the word problem} 

In \cite{grigorchukcrypto} the authors proposed a public key cryptosystem based on the word decision problem in the family of Grigorchuk groups. This protocol is similar to one constructed by Wagner and Magyarik \cite{WaMa85}, and it is the first application of groups acting on rooted trees in cryptography. The protocol works as follows:

Let $G_{\omega}$ be a Grigorchuk group, where $\omega \in \{0,1,2\}^{\mathbb{N}}$.
\begin{itemize}
    \item Alice chooses her private element $\omega \in \{0,1,2\}^{\mathbb{N}}$. The public key of Alice is a finite set of words that denote the identity element of $G_{\omega}$, and two different words $a, b$. 
    \item Bob encrypts the message by choosing a bit $i \in \{0,1\}$ and then he obtains a word $w_i^*$ by applying randomly a series of addition of deletions of relations provided by Alice.
    \item Alice decrypts the message by checking whether $w_i=w_0$ or $w_i=1$.
\end{itemize}

In \cite{grigorchukcrypto} the authors claim that this public key cryptosystem is secure since there is not enough information to determine the key $\omega$. The cryptosystem presented above was first attacked by Hofheinz and Steinwandt in \cite{attack1}. However, their algorithm is mainly based on brute force. Later in 2003, Petrides \cite{grigorchukcripto2} claimed that this protocol is vulnerable since the public key given, gives too much information and allows a possible attacker to easily obtain the private key. 



\subsection{Key agreement scheme based on the simultaneous conjugacy search problem}

Anshel-Anshel-Goldfeld key-agreement \cite{AAG} is well-known and based on the complexity of the simultaneous conjugacy search problem. 
In 2019, Grigorchuk and Grigorev proposed certain classes of automaton groups for the AAG protocol. Concrete examples they proposed are the first Grigorchuk group and the Basilica group.

Let $G$ be a group. and consider public elements $a_1, \dots, a_n, b_1, \dots, b_m \in G$. For an element $s \in G$, we write $\overline{s}=s^{\pm 1}$. The AAG key exchange protocol works as follows:
\begin{itemize}
    \item Alice chooses her private element $a = \overline{a}_{p_1} \dots \overline{a}_{p_s} \in \langle a_1, \dots, a_n \rangle$.
    \item Bob chooses his private element $b = \overline{b}_{q_1} \dots \overline{b}_{q_t} \in \langle b_1, \dots, b_m \rangle$.
    \item Alice sends $b_i^{a}=a^{-1}b_i a$, with $1 \leq i \leq m$.
    \item Bob sends  $a_j^{b}=b^{-1}a_jb$, with $1 \leq j \leq n$.
    \item Alice computes $bab^{-1} =b\bar{a}_{p_1}b^{-1}\dots b\bar{a}_{p_s}b^{-1}$.
    \item Bob computes $a^{-1}ba =a^{-1}\bar{b}_{q_1}a^{-1}\dots a^{-1}\bar{b}_{q_t}a$. 
\end{itemize}
The common key is the commutator $[a,b^{-1}]=a^{-1}(bab^{-1}) = (a^{-1}ba)b^{-1}$ computed by Alice and Bob.

An evesdropper has to find elements $A \in \langle a_1,\dots, a_n\rangle$ and $B \in \langle b_1, \dots, b_m\rangle$ such that $A^{-1}b_iA = a^{-1}b_ia$ with $1 \leq  i \leq m$ and $Ba_jB^{-1} = ba_jb^{-1}$ with $1 \leq  j \leq n$, where $a^{-1}b_ia$ and $ba_jb^{-1}$ are known. Finally, one has to check that $a^{-1}bab^{-1}=A^{-1}BAB^{-1}$.

Furthermore, a possible attacker has to search a solution $A \in \langle a_1, \dots, a_n \rangle$ of the problem $A^{-1}b_iA = a^{-1}b_ia$ with $1 \leq  i \leq m$. The security of this system in general lies on the difficulty of variation of conjugacy search problem in certain groups. For example Ko-Lee in \cite{braid} proposed braid groups and Kahrobaei-Eick \cite{Eickconj,polycyclic} proposed polycyclic groups.

\section{A new cryptosystem metascheme based on generic automaton groups}

In this section, we present a new public key metascheme that can be applied in principle to any automaton group. In a couple of group-based cryptosystems, the platform group must have the word decision problem, which is a tractable problem because the exchanged key is an element of the platform group and therefore we must be able to easily verify whether two words represent the same element in that group. In some classes of automaton groups the word decition problem that is a {\bf PSPACE}-complete problem see \cite{PSPACE}, so beyond the {\bf NP}-completeness.  These groups do not possess in general a normal form. Here, we present a metascheme that bypasses this by making use of the action of an automaton group. These groups act on infinite rooted trees in a quite rich and involuted way, since two equivalent words act in the same way, we may use this action to encode a message or a common shared key.  

\subsection{The protocol}
The public key is the following tuple
\[
\left(\mathcal{A}, (b_1,\ldots, b_{\ell}),\{r_1, \ldots, r_k\} \right)
\]
where $\mathcal{A}=(Q,X, \cdot, \circ)$ can be any automaton group, $ (b_1,\ldots, b_{\ell})$ is an ordered $\ell$-tuple of words $b_i\in (Q\cup Q^{-1})^*$, and $\{r_1,\ldots, r_k\}$ is a set of words $r_i\in  (Q\cup Q^{-1})^{\star}$ representing the identity $r_i= \mathds{1}$ in the associated group $G(\mathcal A)$.
\begin{itemize}
\item [a)] Alice chooses randomly a keyword $A\in (Q\cup Q^{-1})^{\star}$ and in the free monoid $(Q\cup Q^{-1})^*$ she computes the conjugate $\hat{b}_i=Ab_iA^{-1}$ for each $i=1,\ldots, \ell$.  Alice manipulates each $\hat{b}_i$ by rewriting it using the semi-Thue system $\leftrightarrow_R^*$ associated to the set $R=\{r_1,\ldots, r_k\}$ defined in Section~\ref{sec: preliminaries}. Starting from $\hat{b}_i$ Alice randomly applies productions of $\leftrightarrow_R^*$ 
and rewrites $\hat{b}_i$ into a new (equivalent) word $c_i$ with $\hat{b}_i\leftrightarrow_R^* c_i$.  
\item [b)] Alice sends the ordered tuple $(c_1,\ldots, c_{\ell})$ to Bob.
\item [c)] Let $B=\{b_1,\ldots, b_{\ell}\}$. Bob computes randomly a sufficiently large key $u\in (B\cup B^{-1})^*$ such that $A$ does not commute with $u$. 
Since the word problem is in general {\bf PSPACE}-complete, instead of checking if $u^A=u$ in $G(\mathcal A)$, Bob tests for some words $w\in X^*$ whether $u^A\circ w\neq u\circ w$, if it fails this test for several words $w_1,\ldots, w_k$, then he discards $u$ and he repeats the protocol, otherwise $u$ is accepted. Suppose that the generated key is $u=b_{i_1}^{e_1}\ldots b_{i_m}^{e_m}$, for some $e_i\in\{1,-1\}$ and $b_{i_j}\in (B\cup B^{-1})$. Bob computes the corresponding word:
\[
u^A:= c_{i_1}^{e_1} \ldots c_{i_m}^{e_m}
\] 
which is now equal to $A uA^{-1}$ in $G(\mathcal A)$. Finally, he sends $u^A$ to Alice.
\item [d) ]Alice knowing her private key $A$ can consider the word 
\[
U=A^{-1} u^A A
\]
that represents the same element as $u$ in $G(\mathcal{A})$. Now, Alice and Bob may use their words $U, u$, respectively to encrypt any word $m\in X^*$. 
\item [e)] If Bob wants to send the word $m\in X^*$ to Alice, he encodes $m$ by taking $M=u\circ m$ and sends it back to Alice. Then Alice recovers the message by computing 
\[
m=U^{-1} \circ M
\]
Conversely, if Alice wants to send $m\in X^*$ to Bob, she has to encode the message by computing $M=U\circ m$, and she sends it to Bob. Then Bob will calculate $u^{-1}\circ M=m$ to recover it.
\end{itemize}

\subsection{Cryptanalysis}

In this section, we observe and highlight some properties of the automaton groups that could make the protocol secure and also propose some open problems to study the security of such a metascheme in some critical steps where a possible attacker could operate.
\\
First, we note that when Alice computes the conjugates $\hat{b}_i=Ab_iA^{-1}$ in the free monoid $(Q\cup Q^{-1})^{\star}$ for $i=1,\ldots, \ell$, all these elements represent group elements $\beta_i\in G(\mathcal A)$. If Alice would send these words as they are through the public channel, knowing $b_i$ and $\hat{b}_i$ an attacker could easily try to guess the secret keyword $A$ that satisfies the equation $\hat{b}_i=Ab_iA^{-1}$ in the free monoid $(Q\cup Q^{-1})^*$. For this reason, Alice changes each $\hat{b}_i$ by rewriting it using the semi-Thue system $\leftrightarrow_R^*$ associated to the set $R=\{r_1,\ldots, r_k\}$. By randomly rewriting $\hat{b}_i$ into words $c_i$, both $\hat{b}_i$ and $c_i$ represents the same element in $G(\mathcal A)$. However, in general, the set of words $X\in \{\hat{b}_1,\ldots, \hat{b}_{\ell}\}$ such that $X\leftrightarrow_R^n c_i$ that requires $n$ steps to reach $c_i$, may grow exponentially with respect to the parameter $n$, making the problem of recovering the word $\hat{b}_i$ from $c_i$ intractable (as well as recovering the key $A$).
In step c) of our protocol it is important that Bob chooses a word $u$ that does not commute with the key $A$ in $G(\mathcal A)$, this is a critical step since in this case he would send in the public channel the key $u$. As we have already pointed out since the word problem for an automaton group is in general a {\bf PSPACE}-complete (see \cite{PSPACE}), Bob cannot test whether $u^A= u$ in $G(\mathcal A)$ or not. The test performed in step c) ensures to have a key $u$ with the aforementioned property. In the authors' knowledge, there is no theoretical result regarding the probability of randomly choosing a word in an automaton group that does not commute with a fixed element $A$, we strongly believe such probability to be in general strictly positive for non-commutative automaton groups.

\begin{problem}
Is it true that an automaton group $\mathcal{A}$ that defines a non-commutative group $\mathcal{G}(\mathcal{A})$, the probability of choosing uniformly two words $u, A$ of length less or equal to $N$ that does not commute is strictly positive for any $N$ sufficiently large? In case it is not true, characterize such groups. 
\end{problem}
We point out that in general, $G(\mathcal{A})$ does not have a normal form, and the word problem is beyond intractability which makes the previous metascheme suitable to be implemented with platform (automaton) groups which may be quite robust to classical attacks of the cryptosystems based on groups. This is at the cost of Alice and Bob sharing the same group element but with different representatives: $U$ for Alice, and $u$ for Bob. However, our protocol overcomes this issue by making use of the automaton's natural action to encrypt and decrypt a word $m\in X^*$ (this is done in steps d) and e)). 
This is a critical step since we are using the natural action of the automaton to encrypt the word $m$, so we may directly encode the message or a common private key to be used in a possible symmetric protocol. In general, the action of an automaton group can be quite complicated, and obtaining some information of $m$ from $M=u\circ m$ seems very difficult, especially if there is not so much information regarding the secret key $u$. In this direction, an interesting general question regards the search space for the message $m$ knowing the encrypted one $M$. 
\begin{problem}
Given an automaton group $\mathcal{A}$, and a word $M$, give some interesting upper bounds (for instance for the platform group proposed in Section~\ref{sec: platform for our protocol}) to the set 
\[
\{m': v\circ m'=M\mbox{ for all }v\in (Q\cup Q^{-1})^*\mbox{ with length less or equal to n}\}
\]
where $n$ is the length of the word $u^A$ sent by Bob to Alice in step c). The information that a possible attacker has, is the word $u^A$ which is a conjugate of the key $u$ used to encrypt $m$, so it would be interesting to study the average of the Hamming distance between the encrypted message $M$ and the words $h\circ M$ where $h\in \mathcal{G}(\mathcal{A})$ are elements whose inverses $h^{-1}$ are conjugated to $u^A$ in $\mathcal{G}(\mathcal{A})$.
\end{problem}
Another point of strength of our protocol is the fact that any attack on the AAG protocol may be much harder in our context since, as we will see later, the conjugacy decision problem is undecidable. 
A general attack on the AAG protocol is to try to recover the element $\hat{b_i}$ from the representative $c_i$ sent through the channel. In this case, knowing $\hat{b_i}$ and $b_i$ it would be easy to find the word $A$ satisfying the equation $Ab_iA^{-1}= \hat{b_i}$ in the free monoid. However, as we have pointed out before, from the authors' perspective, it appears a hard task to recover a word from another obtained by randomly inserting some known factors. Nevertheless, this issue deserves a deeper analysis, for instance by calculating the dimension of the search space of the set $X\in \{\hat{b}_1,\ldots, \hat{b}_{\ell}\}$ such that $X\leftrightarrow_R^n c_i$ for instance for the platform group described in Subsection~\ref{sec: platform for our protocol}. Furthermore, the classical conjugacy search attack on the Anshel-Anshel-Goldfeld protocol, i.e., trying to solve the equation $Xb_iX^{-1}=c_i$ in the automaton group $G(\mathcal A)$ strongly depends on the decidability of the conjugacy decision problem, see the reference \cite{VS} comparing the conjugacy decision problem and conjugacy search problem.
Since there are automaton groups with undecidable conjugacy decision problem \cite{SUNIC2012148}, we deduce that in general, such an attack would be beyond the intractability, making this kind of attack impossible even for a quantum computer. 

\subsection{Platform groups}\label{sec: platform for our protocol}
We emphasize that our protocol is not based on a specific automaton group, since $\mathcal{A}$ is part of the public key, so from one side it is very difficult to perform an attack for a generic automaton $\mathcal{A}$. On the other side, to assess the security of such protocol, it is important to consider a class of automaton groups that may be used as a platform. In this section, we propose such a class. These groups have been studied in \cite{SUNIC2012148} and it contains an example of an automaton group with undecidable conjugacy problem and for completeness we will describe it below.
Let $n>1$ and consider $G_n \leq \mbox{ Aff}_d(\mathbb{Z}_n)$ the group consisting of some affine invertible transformations over the ring $\mathbb{Z}_n^d$, where $\mathbb{Z}_n$ is the ring of $n$-adic integers. It has been proved that
$$
G_n \cong \langle a_1, a_2, \dots, a_d, t | [a_i, a_j]=1, \ t a_j t^{-1}=a_1^{m_{1,j}}\dots a_d^{m_{d,j}}, \ 1 \leq i, j \leq d \rangle.
$$
The authors proved that $G_n$ can be realized by a finite automaton acting on an alphabet of size $n^d$  and it has in general undecidable conjugacy problem.
\section{Post-quantum Analysis}

Group-based cryptography is a relatively unexplored family in post-quantum cryptography. However, the complexity of the Hidden Subgroup Problem (HSP) in automaton groups and its relationship to more well-known hardness problems, particularly with respect to its security against quantum adversaries, has not been well understood and is a significant open problem for researchers in this area. Kuperberg in \cite{kuperberg2020} has conjectured that ``HSP is probably hard for most infinite groups, but they have a wide variety of behaviors''. Kuperberg made an announcement at \'Ecole Normale Sup\'erieure Un apr\'es-midi de th\'eorie des groupes in Ulm Paris, January 10, 2023 that ``if $G$ is a non-abelian free group and $H$ is assumed to be normal, then the hidden subgroup problem is NP-hard and a Shor-type algorithm is implausible at best. The proof depends on the structure theory of small-cancellation groups.''
The Basilica groups and Brunner-Sidki-Viera automaton groups are both infinite torsion-free classes of automaton groups. For such groups, according to Kuperberg's conjecture, the HSP could be NP-hard.

The first Grigorchuk automaton groups and more generally, of all just infinite groups (i.e. all infinite groups whose proper quotients are finite) such as the classes of Gupta-Sidki $p$-groups, for p odd prime, have the properties that every subgroup is of finite index. For other classes of automaton groups, for example the Basilica automaton groups, all maximal subgroups are of finite index.

Since the existence of finite subgroups is essential to a HSP based attack, these cases should be investigated more carefully.

In addition, some progress has been made in quantum solutions to HSP for certain nonabelian finite groups, such as semidirect products of abelian groups, or groups with the property that all subgroups are normal. Kuperberg in \cite{kuperberg2005subexponential}, and Regev in \cite{regev2004subexponential} give a subexponential-time quantum algorithm for the dihedral hidden subgroup problem. For a slower but less space-expensive algorithm, one can also use a generalised version of an algorithm due to Regev \cite{regev2004subexponential}. The generalised version appears in Theorem~5.2 in \cite{childs2014constructing}. In \cite{battarbee2022subexponential}, the authors give the first dedicated security analysis of Semidirect Discrete Logarithm Problem (SDLP). In particular, they provide a connection between SDLP and group actions, a context in which quantum subexponential algorithms are known to apply. They are therefore able to construct a subexponential quantum algorithm for solving SDLP, thereby classifying the complexity of SDLP and its relation to known computational problems.
\section{Conclusion and Open Problems}

In this paper, we have presented the current state of automaton group-based cryptography.  We conclude by collecting a few open problems below with the hope of stimulating interest in their solutions.


\begin{enumerate}
\item What is the computational complexity of the (simultaneous) conjugacy decision/search problem in other classes of automaton groups?
    \item Are there other cryptographic schemes that could be based on hard problems in automaton groups?
    \item The rigorous cryptographic security analysis of the proposed cryptosystems are of interest for the cryptographic community.
    \item There are various attacks such as Length-based-Attack, quotient attacks, etc have been analysed for braid groups \ \cite{GarberLBA06} and polycyclic groups \cite{garber2013analyzing}. To ensure the security of the proposed cryptosystems based on difficulty of conjugacy search problems in automaton groups, such attacks should be implemented and tested for these groups.
    \item Solving the HSP for different classes of automaton groups, will give special interest for the post-quantum community.
\end{enumerate}

\vspace{0.2cm}

\bibliographystyle{plain}
\bibliography{bib}

\end{document}